\documentclass[letterpaper,english,american,aps, pre, twocolumn, showkeys, showpacs, superscriptaddress]{revtex4-1}
\usepackage[T1]{fontenc}
\usepackage{babel}

\usepackage{units}
\usepackage{amsmath}
\usepackage{graphicx}
\usepackage{amssymb}
\usepackage[unicode=true, 
 bookmarks=false,
 breaklinks=false,pdfborder={0 0 1},backref=false,colorlinks=false]
 {hyperref}
\hypersetup{pdftitle={Non-equilibrium dynamics of stochastic point processes with refractoriness},
 pdfauthor={Deger, Helias, Cardanobile, Atay, Rotter},
 colorlinks=true,linkcolor=black,citecolor=black,urlcolor=black,filecolor=black}
 
\makeatletter

\AtBeginDocument{}

\DeclareFontEncoding{LGR}{}{}

\@ifundefined{textcolor}{}
{%
 \definecolor{BLACK}{gray}{0}
 \definecolor{WHITE}{gray}{1}
 \definecolor{RED}{rgb}{1,0,0}
 \definecolor{GREEN}{rgb}{0,1,0}
 \definecolor{BLUE}{rgb}{0,0,1}
 \definecolor{CYAN}{cmyk}{1,0,0,0}
 \definecolor{MAGENTA}{cmyk}{0,1,0,0}
 \definecolor{YELLOW}{cmyk}{0,0,1,0}
 }

\makeatother

\begin{document}
\global\long\def\defeq{\stackrel{\mathrm{def}}{=}}
\global\long\def\autocorr{R}
\global\long\def\gammarate{\beta}
\global\long\def\expect{\mathbb{E}}

\title{Non-equilibrium dynamics of stochastic point processes with refractoriness}

\author{Moritz \surname{Deger}}

\altaffiliation{M.D. and M.H. contributed equally to this work; Electronic address: deger@bcf.uni-freiburg.de}

\affiliation{Bernstein Center Freiburg \& Faculty of Biology, Albert-Ludwig University,
79104 Freiburg, Germany}

\author{Moritz \surname{Helias}}

\altaffiliation{M.D. and M.H. contributed equally to this work; Electronic address: deger@bcf.uni-freiburg.de}

\affiliation{RIKEN Brain Science Institute, Wako City, Saitama 351-0198, Japan}

\author{Stefano \surname{Cardanobile}}

\affiliation{Bernstein Center Freiburg \& Faculty of Biology, Albert-Ludwig University,
79104 Freiburg, Germany}

\author{Fatihcan M. \surname{Atay}}

\affiliation{Max Planck Institute for Mathematics in the Sciences, 04103 Leipzig,
Germany}

\author{Stefan \surname{Rotter}}

\affiliation{Bernstein Center Freiburg \& Faculty of Biology, Albert-Ludwig University,
79104 Freiburg, Germany}

\date{\today}

\begin{abstract}
Stochastic point processes with refractoriness appear frequently in
the quantitative analysis of physical and biological systems, such
as the generation of action potentials by nerve cells, the release
and reuptake of vesicles at a synapse, and the counting of particles
by detector devices. Here we present an extension of renewal theory
to describe ensembles of point processes with time varying input.
This is made possible by a representation in terms of occupation numbers
of two states: Active and refractory. The dynamics of these occupation
numbers follows a distributed delay differential equation. In particular,
our theory enables us to uncover the effect of refractoriness on the
time-dependent rate of an ensemble of encoding point processes in
response to modulation of the input. We present exact solutions that
demonstrate generic features, such as stochastic transients and oscillations
in the step response as well as resonances, phase jumps and frequency
doubling in the transfer of periodic signals. We show that a large
class of renewal processes can indeed be regarded as special cases
of the model we analyze. Hence our approach represents a widely applicable
framework to define and analyze non-stationary renewal processes.
\end{abstract}
\maketitle

\section{Introduction\label{sec:Introduction}}

\begin{figure}
\begin{center}
\includegraphics[width=0.7\columnwidth]{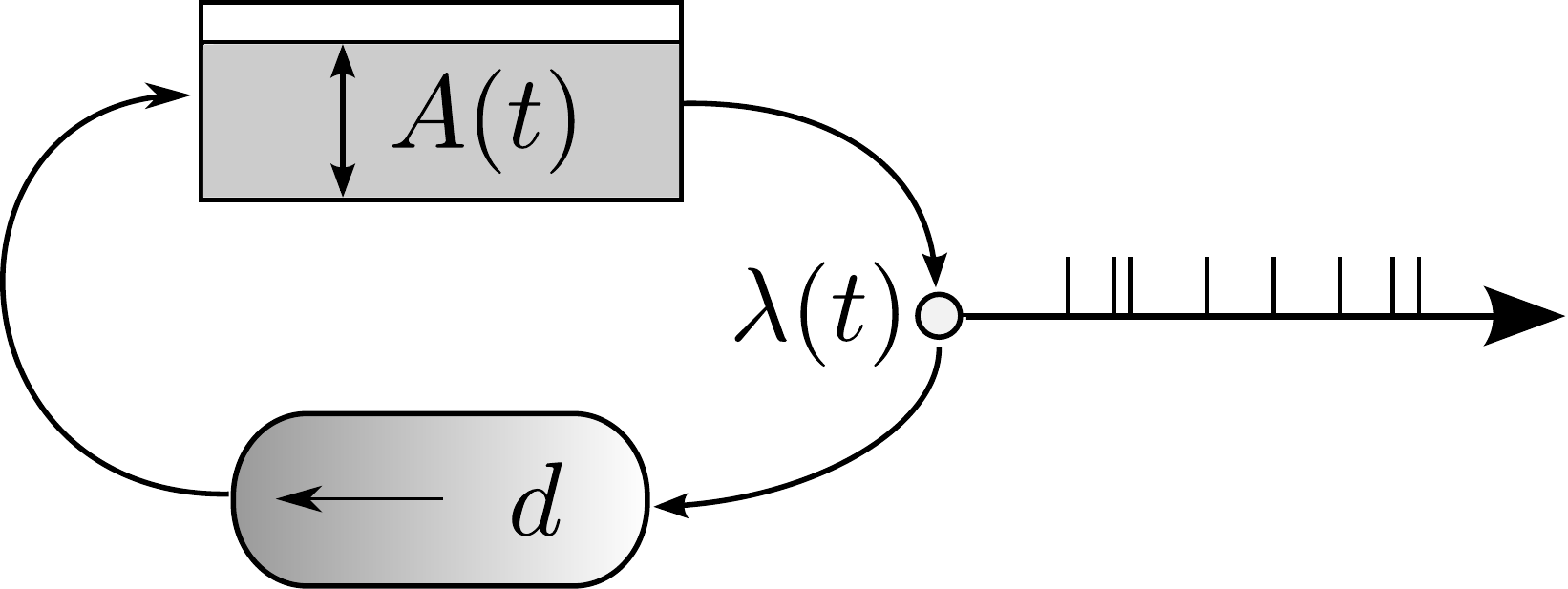}
\end{center}
\caption{Scheme of the ensemble description of the Poisson process with refractoriness: Active component processes produce events with rate $\lambda(t)$ and remain refractory for the duration $d$, illustrated by the delay line. After the dead-time they become active again. The fraction of active component processes is given by $A(t)$.}
\label{fig:ppd_scheme}
\end{figure}

Point processes are stochastic models for time series of discrete
events: a particle passes through an apparatus, a photon hits a detector,
or a neuron emits an action potential \citep{Cox67,LowenTeich05}.
As diverse as these examples are, they share three basic features
that need to enter a statistical description and which are illustrated
in Fig.~\ref{fig:ppd_scheme}. The first feature is refractoriness.
Technical devices to detect point events typically cannot discriminate
events in arbitrarily short succession. This is addressed as the dead-time
of the detector \citep{Grupen08,Muller81}. The process of vesicle
release and transmitter recycling in the synaptic cleft is of similar
nature \citep{Loebel09}. Upon the arrival of an action potential
at the synapse, a vesicle might fuse with the membrane and release
its contents into the synaptic cleft. Subsequently the vesicle is
reassembled for future signaling, but it is available only after a
certain delay, equivalent to a refractory signalling component. In
neurons, refractoriness can be the result of the interplay of many
cellular mechanisms, and possibly also of network effects \citep{Gerstner02}.
In case of cortical neurons, which are driven to produce an action
potential mainly by fluctuations of the input currents \citep{Kuhn04},
refractoriness can model the time it takes to depolarize the membrane
from a hyper-polarized level that follows the action potential into
a range in which action potentials can be initiated by fluctuations.
Generally, refractoriness can be described as a duration $d$ for
which the component cannot be recruited to generate another event.
In Fig.~\ref{fig:ppd_scheme} it is illustrated as a delay line.
After the refractory time is elapsed, the component reenters the pool
of active components that can generate an event. The existence of
such a pool is the second common property of the examples. Each component
process of the ensemble can be either active or refractory. So an
ensemble of neurons, vesicles, or detectors, can be treated in terms
of the occupation of two states, {}``active'' and {}``refractory'',
as depicted in Fig.~\ref{fig:ppd_scheme}, where $A(t)\in[0,1]$
describes the fraction of components which are active at time $t$
and $1-A(t)$ is the fraction of components that are currently refractory.
The third feature is the stochastic nature of event generation. The
time of arrival of a particle at a detector, the fluctuation of the
membrane potential of a neuron that exceeds the threshold for action
potential initiation, and the release of a vesicle into the synaptic
cleft can under many conditions be assumed to happen stochastically.
Given an independent transition density of $\lambda(t)$ per time
interval, event generation follows an inhomogeneous Poisson process,
as indicated in Fig. \ref{fig:ppd_scheme}. In the example of a detector,
$\lambda(t)$ corresponds to the actual rate of incoming particles,
and we will call it the input rate in the following. We distinguish
two models of systems which share the properties described above:
If the refractoriness has a fixed duration we obtain the well-known
Poisson process with dead-time (PPD). If the duration is drawn randomly
from a specified distribution we call the model the Poisson process
with random dead-time (PPRD). 

In the following we describe an extension of renewal theory for ensembles
of point processes with time varying input. For stationary input rate,
many previous publications have investigated the statistics of the
PPD \citep{Brenguier89,Yu00,Teich76,Picinbono09}. In case of slowly
varying input rates, expressions for mean and variance of detector
counts have been derived \citep{Vanucci78}, and recently a method
was proposed to correct for clustered input events \citep{Larsen09}.
PPDs with non-stationary input rate have been employed as a model
for the signal transduction in the auditory nerve \citep{Turcott94}.
A sudden switch of the input rate was found to induce strong transients
of the ensemble output rate. These reflect a transiently perturbed
equilibrium of the occupation numbers in our two-state model, and
we quantitatively analyze this case for the PPD. Response transients
to rapidly changing input, in fact, explain the relation between neural
refractoriness and neural precision \citep{Berry98}. Furthermore,
periodic input profiles are known to be distorted by refractoriness
\citep{Johnson83}. Here we derive the mapping of periodic input to
output in the steady state and uncover the impact of refractoriness
on the transmission. Interacting populations of refractory neurons
have been studied in \citep{Wilson72_1}. This approach, in contrast
to ours, neglects the effects of refractoriness on short time scales
due to temporal coarse-graining of the population dynamics. 

From a more abstract perspective, the PPD is a very simple example
of a point process that exhibits stochastic transients, which are
not shared by the ordinary Poisson process. Besides its many applications,
the PPD therefore is a prototype system to study non-equilibrium phenomena
in point process dynamics. Generally, non-stationary point processes
can be defined by two different models: by rescaling time \citep{Brown01,Reich98,Nawrot08_374}
or by time-dependent parameters of the hazard. The drawback of the
former method is that the transformation from operational to real
time distorts the inter-event intervals, such that, for example, a
constant refractory period is not maintained. An example how a time-dependent
hazard function can be derived from an underlying neuron model with
time-dependent input is given in \citep{Muller07_2958}. Our approach
differs with regard to the choice of the hazard function, which enables
rigorous analysis of the dynamics of the process.

To analytically investigate non-equilibrium phenomena in ensembles
of renewal process, a typical approach is to use a partial differential
equation (PDE) for the probability density of the ages of the components
(time since the last event)\citep{Gerstner02}. In Section \ref{sec:dynamics_ppd}
we derive the two-state representation of the PPD from the dynamics
of the age density. We present analytical solutions of the population
dynamics for the response to a step change in the input rate in Section
\ref{sec:solution_ppd_step}, and to periodic input rate profiles
in Section \ref{sec:periodic_input_ppd}. Finally, in Section \ref{sec:PPRD},
we generalize our results to random refractoriness. We compute the
effective hazard function of the resulting inhomogeneous renewal process,
connecting it to the framework of renewal theory. For the PPRD with
gamma-distributed dead-times, as applied recently to model neural
activity \citep{Toyoizumi09}, we show how the dynamics in terms of
a distributed delay differential equation can be reduced to a system
of ordinary (non-delay) differential equations. Again we study the
transient response of an ensemble of processes to a step-like change
in the input rate and the transmission of periodic input. We observe
that both distributed and fixed refractoriness lead to qualitatively
similar dynamical properties. At last we identify the class of renewal
processes that can be represented as a PPRD. As it turns out, this
covers a wide range of renewal processes.

\section{Dynamics of an ensemble of PPDs\label{sec:dynamics_ppd}}

Point processes can be defined by a hazard function \begin{equation}
h(t,\mathcal{H}_{t})\defeq\lim_{\epsilon\to0}\frac{1}{\epsilon}\mathbb{P}\left[\text{event in }[t,t+\epsilon]\,|\,\mathcal{H}_{t}\right],\label{eq:hazard_general}\end{equation}
which is the conditional rate of the process to generate an event
at time $t$, given the history of event times $\mathcal{H}_{t}$
up until $t$. A process is a renewal process \citep{Cox67} if the
hazard function depends only on the time $\tau$ since the last event
(age) instead of the whole history $\mathcal{H}_{t}$. This can be
generalized to the inhomogeneous renewal process which, additionally,
allows for an explicit time dependence of the hazard function $h(t,\mathcal{H}_{t})=h(t,\tau)$.

Here we consider an ensemble of point processes defined by the hazard
function \begin{equation}
h(t,\tau)=\lambda(t)\theta(\tau-d),\label{eq:def_hazard}\end{equation}
where $\theta(t)=\{1\text{ for }t\geq0,\;0\text{ else\}}$ denotes
the Heaviside function, $d\geq0$ is called dead-time, $\lambda(t)\geq0$
is the time-dependent input rate and $\tau\geq0$ is the age of the
component process. This is an inhomogeneous renewal process, which
is known as the Poisson process with dead-time (PPD). The state of
an ensemble of such processes can be described by the time-dependent
probability density of ages $a(t,\tau)$, for which a partial differential
equation is known \citep{Gerstner02} \begin{equation}
\frac{\partial}{\partial t}a(t,\tau)=-\frac{\partial}{\partial\tau}a(t,\tau)-h(t,\tau)a(t,\tau).\label{eq:PDE}\end{equation}
Solutions must conserve probability, which manifests itself in the
boundary condition $a(t,0)=\nu(t)$, with the event rate of the ensemble
\begin{equation}
\nu(t)\defeq\int_{0}^{\infty}h(t,\tau)a(t,\tau)\, d\tau=\lambda(t)A(t).\label{eq:def_rate}\end{equation}
In the second step we inserted \eqref{eq:def_hazard} and introduced
the active fraction of component processes with age $\tau\geq d$
\begin{equation}
A(t)\defeq\int_{d}^{\infty}a(t,\tau)\, d\tau.\label{eq:def_A}\end{equation}
For $\tau<d$, Eq.~\eqref{eq:PDE} simplifies to $\partial_{t}a(t,\tau)=-\partial_{\tau}a(t,\tau)$,
implying \begin{equation}
a(t+u,u)=a(t,0)=\nu(t)\quad\forall u\in[0,d).\label{eq:delay_prop}\end{equation}
Since $a(t,\tau)$ is normalized we obtain with the boundary condition
and \eqref{eq:delay_prop}, \begin{equation}
1=\int_{0}^{\infty}a(t,\tau)\, d\tau=\int_{t-d}^{t}\nu(s)\, ds+A(t).\label{eq:normalization_A}\end{equation}
This equation is the starting point of the analysis of interacting
populations of refractory neurons in \citep{Wilson72_1}. Differentiation
of \eqref{eq:normalization_A} by $t$ yields \begin{equation}
\frac{d}{dt}A(t)=\lambda(t-d)A(t-d)-\lambda(t)A(t)\;,\label{eq:diffeq_A}\end{equation}
which is a linear delay differential equation (DDE) with time-dependent
coefficients. Its forward solution for input $\lambda(t)$ is uniquely
defined given $A(t)$ on an interval of length $d$ \citep{Bellman63}.
However, not all solutions of \eqref{eq:diffeq_A} can be interpreted
physically, since by differentiation of \eqref{eq:normalization_A}
additive constants are lost. Only if the initial trajectory satisfies
\eqref{eq:normalization_A}, Eq.~\eqref{eq:diffeq_A} determines
the time evolution of the ensemble. With \eqref{eq:def_rate}, the
time-dependent output rate $\nu(t)$ follows. Note that only in the
case of the {}``pure'' Poisson process with $d=0$ we obtain $\nu(t)=\lambda(t)$,
because $A(t)=1$ by \eqref{eq:normalization_A}.  

Eq.~\eqref{eq:diffeq_A} represents a more accessible description
of the process in terms of the occupation of the active and the refractory
state (see Fig.~\ref{fig:ppd_scheme} and Section \ref{sec:Introduction})
compared to the dynamics of the probability density of ages \eqref{eq:PDE}.
This description is feasible because of the particular nature of the
hazard function of the PPD \eqref{eq:def_hazard}. In the following
we will consider specific solutions of \eqref{eq:diffeq_A}.

\section{Solutions for a step input\label{sec:solution_ppd_step}}

If $\lambda(t)=\lambda$ is constant, given the occupation $A(t)=u(t)$
on the first interval $t\in[-d,0]$ with $u:[-d,0]\to[0,1]$, solutions
to \eqref{eq:diffeq_A} are known in integral form \citep{Bellman63}\begin{equation}
A(t)=u(0)g(t)+\int_{0}^{d}\lambda u(s-d)\, g(t-s)\, ds\label{eq:solution_initial_u}\end{equation}
for $t\geq0$, where we introduced the fundamental solution $g(t)$.
It obeys $g(t)=\{0\text{ for }t<0,\;1\text{ for }t=0\}$ and solves
\eqref{eq:diffeq_A} for $t>0$. As we will show, here $g$ is in
fact the shifted and scaled auto-correlation function $\autocorr$
of the process. 

The inter-event interval density of the stationary PPD is $f(t)=\lambda\theta(t-d)e^{-\lambda(t-d)}$.
For $t\geq d$, the integral equation\begin{equation}
A(t)=(f\star A)(t),\label{eq:conv_invariance}\end{equation}
is equivalent to the delay differential equation \eqref{eq:diffeq_A},
which can be proven by differentiation with respect to $t$ ($\star$
denotes the convolution). The auto-correlation function \citep{Cox67}
\begin{equation}
\autocorr(t)=\sum_{k=0}^{\infty}f^{\star k}(t),\label{eq:autocorr}\end{equation}
with $f^{\star k}(t)\defeq(f^{\star(k-1)}\star f)(t)$ for $k\geq1$
and $f^{\star0}\defeq\delta(t)$, solves \eqref{eq:conv_invariance}
for $t\geq d$, and hence is a solution of \eqref{eq:diffeq_A} in
that domain. We find for $k\geq1$ that \begin{equation}
f^{*k}(t)=\lambda^{k}(t-kd)^{k-1}e^{-\lambda(t-kd)}\theta(t-kd)/(k-1)!\,.\label{eq:f_k-fold-conv}\end{equation}
Given the initial trajectory $g(t)$ for $t\leq0$, solving \eqref{eq:diffeq_A}
by variation of constants for $t\in[0,d]$ yields $g(t)=\lambda^{-1}f(t+d)=\lambda^{-1}\autocorr(t+d)$.
Then due to uniqueness of the solution it holds for all $t\geq0$
that \begin{equation}
g(t)=\lambda^{-1}\autocorr(t+d)\;.\label{eq:autocorr_fundamental_sol}\end{equation}

\begin{figure*}
\begin{center}
\begin{tabular}{l l}
A & B\\
\includegraphics[scale=1.1]{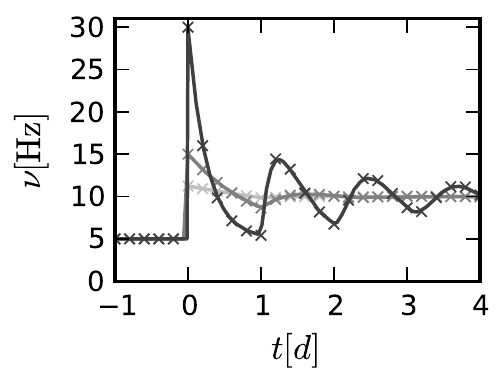} & \includegraphics[scale=1.1]{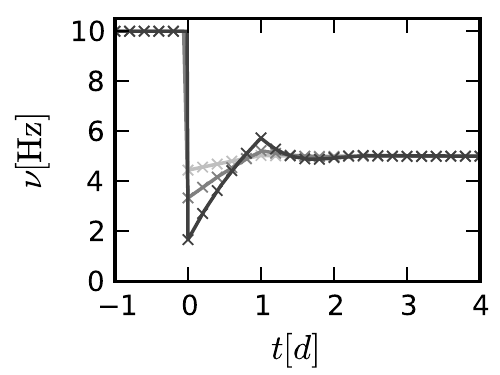}
\end{tabular}
\end{center}
\caption{Transients upon step change of the input rate $\lambda(t)$ at $t=0$. Exact analytical result \eqref{eq:solution_rate_step} (solid lines) and simulation of the ensemble rate of $10^{10}$ processes (crosses). Parameters: $d\,[\mathrm{s}]:$ $0.02,\,0.05,\,0.08$ (light gray, mid gray, dark gray) A: $\lambda_{0}=((5\mathrm{Hz})^{-1}-d)^{-1},\lambda=((10\mathrm{Hz})^{-1}-d)^{-1}$, B: $\lambda_{0}=((10\mathrm{Hz})^{-1}-d)^{-1},\lambda=((5\mathrm{Hz})^{-1}-d)^{-1}$.}
\label{fig:step_transients}
\end{figure*}

We apply these results to compute the response of $A(t)$ if the input
rate is switched from $\lambda_{0}$ to $\lambda$ at $t=0$, given
the process was in equilibrium for $t\le0$. Eq.~\eqref{eq:normalization_A}
determines this equilibrium to  $A(t)\equiv a_{0}=(1+\lambda_{0}d)^{-1}$,
$t\leq0$. In this case, the step change in $\lambda(t)$ enters \eqref{eq:solution_initial_u}
as \begin{equation}
A_{\mathrm{step}}(t)=u(0)g(t)+\int_{0}^{d}\underbrace{\lambda(s-d)}_{\lambda_{0}}u(s-d)\, g(t-s)\, ds,\label{eq:solution_initial_u_step}\end{equation}
for $t\geq0$. We insert \eqref{eq:autocorr_fundamental_sol} to obtain
\begin{eqnarray}
A_{\mathrm{step}}(t) & = & \frac{a_{0}}{\lambda}\autocorr(t+d)+\frac{\lambda_{0}a_{0}}{\lambda}\int_{0}^{d}\autocorr(t+d-s)\, ds\nonumber \\
 & = & \frac{a_{0}}{\lambda}\autocorr(t+d)+\frac{\lambda_{0}a_{0}}{\lambda}\left(1-\frac{1}{\lambda}\autocorr(t+d)\right)\nonumber \\
 & = & \frac{a_{0}\lambda_{0}}{\lambda}\left(1+(\lambda_{0}^{-1}-\lambda^{-1})\,\autocorr(t+d)\right),\label{eq:solution_rate_step}\end{eqnarray}
where we used \eqref{eq:normalization_A}, which holds for $g(t)=\lambda^{-1}\autocorr(t+d)$.
Fig.~\ref{fig:step_transients} shows this analytical solution compared
to direct numerical simulation of an ensemble of PPDs upon a step
change of the input rate $\lambda(t)$ at $t=0$. The output rate
displays a marked transient, which increases with the dead-time $d$
and exhibits oscillations of frequency $1/d$.

\section{Transmission of periodic input\label{sec:periodic_input_ppd}}

We now investigate an ensemble of PPDs with an input rate $\lambda(t)\in\mathbb{R}$
that is periodic. If $T$ is its period, we obtain the Fourier series
$\lambda(t)=\sum_{k=-\infty}^{\infty}\Lambda_{k}e^{ik\omega t}$,
with $\omega=\frac{2\pi}{T}$ and $\Lambda_{k}\in\mathbb{C}$. Then
the steady state solution for the active fraction $A(t)$ of the PPD
is also periodic in $T$, so it can be expressed as $A(t)=\sum_{k=-\infty}^{\infty}\alpha_{k}e^{ik\omega t}$
with $\alpha_{k}\in\mathbb{C}$. Inserted into \eqref{eq:normalization_A}
we obtain \begin{equation}
1=\sum_{k,l=-\infty}^{\infty}\Lambda_{l}\alpha_{k}q_{l+k}e^{i(l+k)\omega t}+\sum_{k=-\infty}^{\infty}\alpha_{k}e^{ik\omega t}\label{eq:norm_A_per}\end{equation}
where for $k\neq0$ \[
\int_{t-d}^{t}e^{ik\omega t}\, dt=\frac{1-e^{-ik\omega d}}{ik\omega}e^{ik\omega t}\defeq q_{k}e^{ik\omega t},\]
and $q_{0}\defeq d$. Since the Fourier basis functions $\{e^{ik\omega t},k\in\mathbb{Z}\}$
are mutually orthogonal, we can separate \eqref{eq:norm_A_per} for
different $k$. This yields the infinite dimensional linear system
of equations \begin{equation}
\delta_{k,0}=q_{k}\sum_{l=-\infty}^{\infty}\Lambda_{l}\alpha_{k-l}+\alpha_{k},\; k\in\mathbb{Z}.\label{eq:alpha_k_system_per}\end{equation}
The ensemble averaged output rate of the PPD defined in \eqref{eq:def_rate}
 then follows as $\nu(t)=\lambda(t)A(t)=\sum_{k=-\infty}^{\infty}\beta_{k}e^{ik\omega t}$
with the spectrum \begin{equation}
\beta_{k}=\sum_{l=-\infty}^{\infty}\alpha_{k-l}\Lambda_{l}=q_{k}^{-1}\left(\delta_{k,0}-\alpha_{k}\right),\label{eq:beta_k(alpha_k,q_k)}\end{equation}
where we used \eqref{eq:alpha_k_system_per}. This given, we replace
the $\alpha_{k}$ by $\beta_{k}$ in \eqref{eq:alpha_k_system_per}
to obtain \begin{equation}
\beta_{k}=\Lambda_{k}-\sum_{l=-\infty}^{\infty}\Lambda_{k-l}q_{l}\beta_{l},\; k\in\mathbb{Z}.\label{eq:beta_k_system_per}\end{equation}
This relation shows how different frequencies of the output rate are
coupled by a convolution with the input spectrum.  Note that inverting
\eqref{eq:beta_k_system_per} yields the spectrum of the time-dependent
input rate $\lambda(t)$ given the spectrum of the output rate signal
$\nu(t)$. 

\begin{figure*}
\begin{center}
\begin{tabular}{l l}
A & B\\
\includegraphics[scale=1.1]{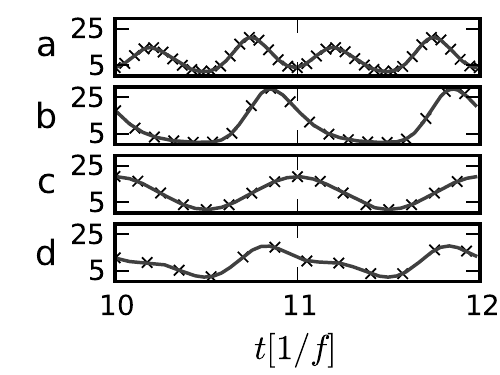} & \includegraphics[scale=1.1]{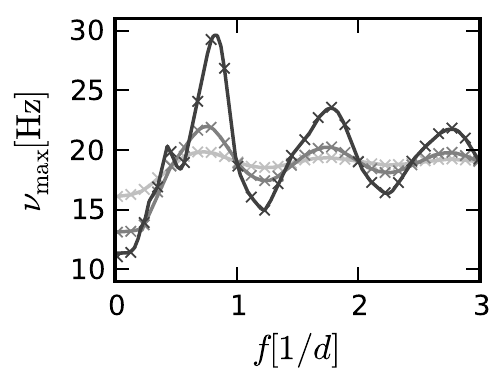}\\[-0.5cm]
C & D\\
\includegraphics[scale=1.1]{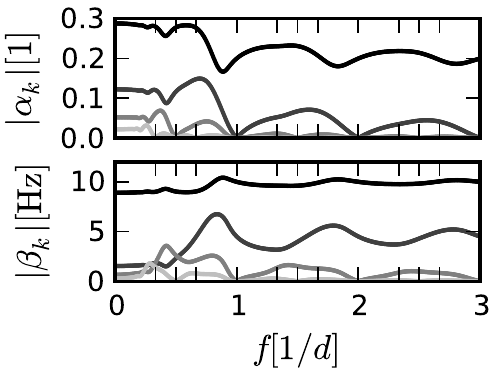} & \includegraphics[scale=1.1]{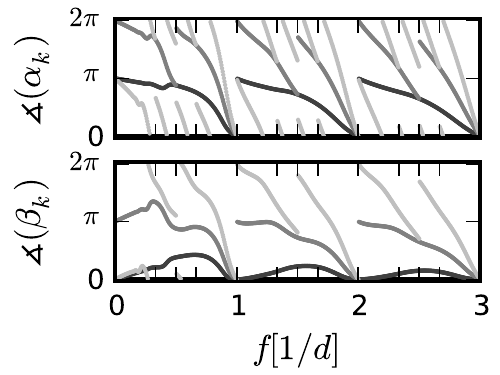} 
\end{tabular}
\end{center}
\caption{Transmission of cosine-modulated input. A, B: Theoretical result (solid lines) and simulation of an ensemble of $10^{10}$ processes (crosses). A: Steady-state rate $\nu(t)$ for different modulation frequencies $f$, with $f d$: $0.42,\,0.85,\,1.0,\,1.4$ (a,b,c,d). B: $\nu_{\max} = \max(\nu(t))$ for different $f$. Here  $d\,[\mathrm{s}]$: $0.02,\,0.05,\,0.08$ (light gray, mid gray, dark gray) C,D: Amplitude (C) and phase (D) of harmonics $k\in\{0,\ldots,3\}$ of $A(t)$ \eqref{eq:beta_k_system_per} (top) and $\nu(t)$ \eqref{eq:beta_k(alpha_k,q_k)} (bottom) as a function of modulation frequency $f$. Grayscale denotes order of harmonics $k:$ $0,1,2,3$ (black, dark gray, mid gray, light gray), $d=80\mathrm{ms}$. Other parameters in A-D: $\lambda(t)=\lambda_{0}(1+0.9\cos(2\pi ft))$, $\lambda_{0}=(\nu_{0}^{-1}-d)^{-1}$, $\nu_{0}=10\mathrm{Hz}$.}
\label{fig:sin_modulation} 
\end{figure*}

Let us now consider the special case of a cosine-modulated input \[
\lambda(t)=\lambda_{0}+\epsilon\cos(\omega t),\]
which we obtain with \foreignlanguage{english}{$\Lambda_{k}=\{0\text{ for }|k|>1\text{, }\frac{\epsilon}{2}\text{ for }k\in\{1,-1\},\,\lambda_{0}\text{ for }k=0\}$},
$\lambda_{0}\geq\epsilon\geq0$. Then for $k\in\mathbb{N}$, \eqref{eq:alpha_k_system_per}
becomes a so-called three-term recurrence relation \citep{Gautschi67}
of the form $0=\alpha_{n+1}+x_{n}\alpha_{n}+y_{n}\alpha_{n-1}$ with
$x_{n}=(q_{n}^{-1}+\lambda_{0})(2/\epsilon)$ and $y_{n}=1$. This
relation has two linearly independent solutions. The unique minimal
solution is convergent and can be obtained from the continued fraction
$r_{n-1}=-y_{n}/(x_{n}+r_{n})$ in a robust manner \citep{Gautschi67}
using the relation $r_{n}=\alpha_{n+1}/\alpha_{n},\, n\geq0$: Setting
$r_{N}=0$ for some $N\in\mathbb{N}$ one computes $(r_{n})_{0\leq n<N}$
backwards and increases $N$ until $r_{0}$ does not change within
the required tolerance. Inserting $\alpha_{1}=r_{0}\alpha_{0}$ into
\eqref{eq:alpha_k_system_per} for $k=0$ we solve for $\alpha_{0}$
to obtain \foreignlanguage{english}{$\alpha_{0}=(1+d(\lambda_{0}+\epsilon\Re(r_{0})))^{-1}$}
(here $\Re$ denotes the real part). The remaining $\alpha_{k}$ follow
recursively from $\alpha_{k+1}=\alpha_{k}r_{k}$ and $\alpha_{-k}=\alpha_{k}^{\star}$,
since $A(t)\in\mathbb{R}$. The spectrum of the output rate is then
given by \eqref{eq:beta_k(alpha_k,q_k)}. Fig.~\ref{fig:sin_modulation}A
shows the output rate $\nu(t)$ for different input rate modulation
frequencies $f=\omega/(2\pi)$. Fig.~\ref{fig:sin_modulation}C,D
display the amplitude and phase of the three lowest harmonics of the
output rate $\nu(t)$ as a function of $f$. The time averaged emission
rate ($\beta_{0}$) depends on the modulation frequency. It is maximized
slightly below the characteristic frequencies $f=k/d$. This is due
to the oscillation of $A(t)$, which is almost in phase at these frequencies
and hence cooperates with the oscillatory hazard rate $\lambda(t)$
to enhance the emission (see Fig.~\ref{fig:sin_modulation}D). Interestingly,
the first ($\beta_{1}$) and second ($\beta_{2}$) harmonic of $\nu(t)$
display maxima at different $f$. At a particular modulation frequency
$f\simeq1/(2d)$ the amplitude of the second harmonic ($\beta_{2}$)
is larger than the first harmonic ($\beta_{1}$), so that the ensemble
activity is effectively modulated with twice the input frequency (see
Fig.~\ref{fig:sin_modulation}A (a) and Fig.~\ref{fig:sin_modulation}C):
the ensemble performs a frequency doubling. Fig.~\ref{fig:sin_modulation}B
shows the maximum over one period of the output rate trajectory. These
maxima are dominated by the maxima of the amplitude of the first harmonic.
In particular, low frequency input signals are transmitted to the
output with strong distortion and reduced intensity, because the fraction
of non-refractory processes, $A(t)$, is in anti-phase (Fig.~\ref{fig:sin_modulation}D)
to $\lambda(t)$ and hence suppresses the output rate's modulation.
This is in contrast to the common view that the PPD transmits slow
signals more reliably than the Poisson process \citep{Gerstner02}.
Note that only if the driving frequency $f=n/d$, $n\in N$ is an
integer multiple of the inverse dead-time then $A(t)=(1+\lambda_{0}d)^{-1}$
is constant in time and the output rate is proportional to $\lambda(t)$
without any distortion.

\section{Random dead-time\label{sec:PPRD}}

For detector devices as well as for neurons, a fixed dead-time might
be a somewhat restricted model. Here we consider the PPRD as described
in the introduction. Upon generation of each event, the PPRD draws
an independent and identically distributed random dead-time with the
probability density function (PDF) $\rho$ for the duration of which
it remains silent. The PPRD is still a renewal process, since it has
no further dependencies on the event history beyond the time since
the last event. As in the case of a fixed dead-time in Section \ref{sec:dynamics_ppd},
the following analysis of the PPRD is based on the conservation of
the total number of processes in an ensemble. Inactive components
must have generated an event at some time in the past, which leads
to the normalization condition \begin{multline}
1=A(t)+\int_{-\infty}^{t}A(t^{\prime})\lambda(t^{\prime})\int_{t-t^{\prime}}^{\infty}\rho(x)dx\, dt^{\prime}\,.\label{eq:normalization_PPRD}\end{multline}
This equation can be seen as the generalization of the normalization
condition \eqref{eq:normalization_A}, from which the DDE \eqref{eq:diffeq_A}
follows by differentiation, to the case of random dead-time. Analogously
the distributed DDE \begin{equation}
\frac{d}{dt}A(t)=-\lambda(t)A(t)+\int_{0}^{\infty}\rho(x)\lambda(t-x)A(t-x)\, dx\label{eq:DDE_pprd}\end{equation}
follows from \eqref{eq:normalization_PPRD} by differentiation with
respect to $t$. Eq.~\eqref{eq:DDE_pprd} describes the time-evolution
of the occupation of the active state for an ensemble of general PPRDs.
Obviously, the dynamics of the PPD \eqref{eq:diffeq_A} is recovered
from \eqref{eq:DDE_pprd} in case of the localized density $\rho(x)=\delta(x-d)$.
In the rest of this section, we will derive the hazard function $h(t,\tau)$
of the PPRD, consider the case of gamma-distributed dead time and
the associated step response, generalize the transmission of periodic
input to random dead-time, and finally identify the class of renewal
processes that can be represented by the PPRD. 

\begin{figure*}
\begin{center}
\begin{tabular}{l l}
A & B\\
\includegraphics[scale=1.1]{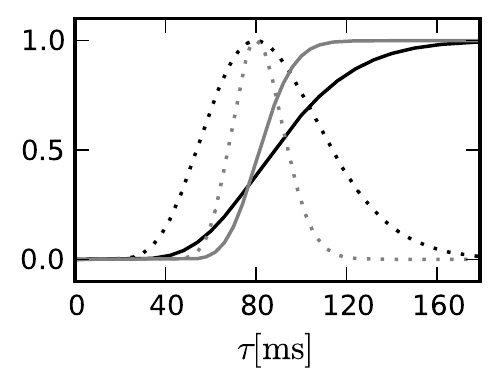} & \includegraphics[scale=1.1]{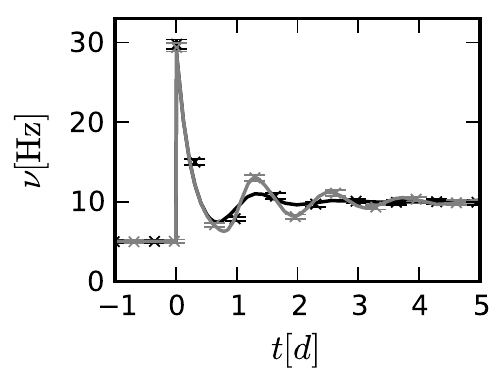}
\end{tabular}
\end{center}
\caption{PPD with random dead-time with mean $80\mathrm{ms}$, shape parameter $n$: $10,\,50$ (black, gray). A: density of dead-times $\rho(\tau)/\max(\rho(\tau))$ \eqref{eq:rho_gamma} (dotted lines) and hazard function $h(\tau)/\lambda_{0}$ \eqref{eq:hazard_pprd} for $\lambda(t)=\lambda_{0}$ (solid lines). B: Transients upon step change of the input rate $\lambda(t)$ at $t=0$. Theoretical result from \eqref{eq:pprd_gamma_sol} (solid lines) and simulation of an ensemble of $10^6$  processes with hazard function \eqref{eq:hazard_pprd} (crosses) averaged over 225 trials. The error bars denote the standard deviation over trials.  $\lambda_{0}=((5\mathrm{Hz})^{-1}-d)^{-1},\,\lambda=((10\mathrm{Hz})^{-1}-d)^{-1}$.}
\label{fig:gen_rec_step}
\end{figure*}

For a given density $\rho(x)$ of the dead-time it is not obvious
what the hazard function of the PPRD is. In order to relate the PPRD
to renewal theory we compute its time-dependent hazard function \eqref{eq:hazard_general}
here. Let \[
Q(t,\tau)\defeq\expect\left[\theta(\tau-x)|\,\text{last event at }t-\tau\right]\]
denote the probability of the process to be active at time $t$, given
the last event occured at $t-\tau$, where $x$ is the random dead-time
and $\expect$ denotes the expectation value with respect to $x$.
The hazard function is then  $h(t,\tau)=\lambda(t)Q(t,\tau)$. With
\begin{multline*}
Q(t,\tau)=\mathbb{P}\left[x<\tau\,|\,\text{last event at }t-\tau\right]\\
=1-\mathbb{P}\left[x\geq\tau\,|\,\text{ev. at }t-\tau\,\cap\,\text{no ev. in }(t-\tau,t)\right]\\
=1-\mathbb{P}\left[x\geq\tau\right]/\mathbb{P}\left[\text{no ev. in }(t-\tau,t)\,|\,\text{ev. at }t-\tau\right]\end{multline*}
we obtain \begin{eqnarray}
h(t,\tau) & = & \lambda(t)\left(1-\mathcal{F}(\tau)/\expect\left[F(t,\tau|x)\right]\right),\label{eq:hazard_pprd}\end{eqnarray}
where \[
\mathcal{F}(\tau)=\int_{\tau}^{\infty}\rho(x)dx\]
is the survivor function of the dead-time distribution and \[
F(t,\tau|x)=\exp(-\theta(\tau-x)\int_{t-\tau+x}^{t}\lambda(t^{\prime})dt^{\prime})\]
is the survivor function of a PPD with dead-time $x$. In case of
constant $\lambda(t)=\lambda$ we further have \[
\expect\left[F(t,\tau|x)\right]=e^{-\lambda\tau}\int_{0}^{\tau}e^{\lambda x}\rho(x)dx+\mathcal{F}(\tau).\]
The hazard function \eqref{eq:hazard_pprd} is shown for constant
$\lambda(t)$ in Fig.~\ref{fig:gen_rec_step}A for the special case
described below. Eq.~\eqref{eq:hazard_pprd} was applied to generate
realizations of the PPRD for Fig.~\ref{fig:gen_rec_step}B.

For gamma-distributed dead-time \eqref{eq:diffeq_A} can be transformed
into a system of ordinary differential equations. We exemplify the
application of \eqref{eq:DDE_pprd} for gamma distributed dead-times
with parameters $n\in\mathbb{N}$ and $\gammarate\in\mathbb{R}^{+}$,
\begin{eqnarray}
\rho(x) & = & \kappa_{n}(x),\label{eq:rho_gamma}\\
\kappa_{n}(x) & = & \gammarate^{n+1}x^{n}e^{-\gammarate x}/n!\,,\label{eq:def_kappa}\end{eqnarray}
with  $\expect[x]=(n+1)/\gammarate$. The time course of the rate
can be obtained from \eqref{eq:DDE_pprd}. Introducing \[
b_{k}(t)\defeq\int_{-\infty}^{t}\kappa_{k}(t-x)\nu(x)\, dx\]
for $0\leq k\leq n$ and $b_{n+1}(t)\defeq A(t)$ and exploiting the
relation \begin{eqnarray*}
\frac{d}{dx}\kappa_{k}(x) & = & \theta(k-1)\gammarate\kappa_{k-1}(x)-\gammarate\kappa_{k}(x)\end{eqnarray*}
for $0\leq k\leq n$ enables to replace the integral in \eqref{eq:DDE_pprd}
by a closed system of ordinary differential equations\begin{eqnarray}
\frac{d}{dt}b_{k}(t) & = & \begin{cases}
-b_{n+1}(t)\lambda(t)+b_{n}(t) & k=n+1\\
\gammarate b_{k-1}(t)-\gammarate b_{k}(t) & 1\leq k\leq n\\
\gammarate b_{n+1}(t)\lambda(t)-\gammarate b_{0}(t) & k=0.\end{cases}\label{eq:pprd_gamma_diffeq_b_k}\end{eqnarray}
For constant $\lambda(t)=\lambda$ this can be written as $\frac{d}{dt}\vec{b}(t)=\mathbf{M}_{\lambda}\,\vec{b}(t)$,
$\vec{b}(t)\in\mathbb{R}^{n+2}$. Hence given the initial state $\vec{b}(0)$
the solution unfolds to \begin{eqnarray}
\vec{b}(t) & = & \exp(\mathbf{M}_{\lambda}t)\cdot\vec{b}(0)\,.\label{eq:pprd_gamma_sol}\end{eqnarray}
With $\nu(t)=\nu=(\lambda^{-1}+\expect[x])^{-1}$ the equilibrium
state follows: Setting the temporal derivatives to $0$ in \eqref{eq:pprd_gamma_diffeq_b_k}
yields $b_{k}=\nu$ for $0\leq k\leq n$, and $b_{n+1}=A=1-\nu\expect[x]$.
The rate response to a switch from $\lambda_{0}$ to $\lambda$ at
$t=0$ is thus given by \eqref{eq:pprd_gamma_sol} where $\vec{b}(0)$
is the equilibrium state for $\lambda_{0}$. A numerical simulation
of the process with gamma distributed refractoriness \eqref{eq:rho_gamma}
with hazard function \eqref{eq:hazard_pprd} and the corresponding
analytical solution \eqref{eq:pprd_gamma_sol} upon a step change
of $\lambda(t)$ are shown in Fig.~\ref{fig:gen_rec_step}. The
simulation of the process was done via rejection \citep{Ripley87}
and averaged over independent runs. The spread of dead-times (Fig.~\ref{fig:gen_rec_step}A)
does not qualitatively change the shape of the response transient
(Fig.~\ref{fig:gen_rec_step}B). 

Analogous to Section \ref{sec:periodic_input_ppd}, we consider the
case of periodic input. We insert the Fourier series of $\lambda(t)$
and $A(t)$ into \eqref{eq:normalization_PPRD} and obtain the same
relation of their spectra \eqref{eq:norm_A_per} as for a single dead-time
 with the altered coefficients \begin{equation}
q_{k}=\int_{0}^{\infty}e^{-ik\omega y}\int_{y}^{\infty}\rho(x)dx\, dy.\label{eq:q_k_PPRD}\end{equation}
As is easily seen, by inserting the localized dead-time PDF $\rho(x)=\delta(x-d)$
the original $q_{k}$ are recovered. Hence all results of Section
\ref{sec:periodic_input_ppd} also hold for the PPRD, but with the
coefficients \eqref{eq:q_k_PPRD}. In particular we would like to
emphasize the validity of the general input-output mapping \eqref{eq:beta_k_system_per}
for arbitrarily distributed dead-time.

Let us now investigate which class of renewal processes can be represented
by the PPRD. We start with an arbitrary renewal process with inter-event
interval $I\in\mathbb{R}_{+}$, defined by its PDF $\iota(x)$. Let
$E\geq0$ be an independent, exponentially distributed interval with
PDF $\epsilon(x)=\lambda e^{-\lambda x}$, and let $R$ be the random
dead-time with PDF $\rho(x)$. For $I$ to be a realization of a PPRD
it must hold for some $\rho$ and $\lambda\geq0$ that  \begin{multline}
I=R+E\Rightarrow\iota=\rho\star\epsilon\Rightarrow\hat{\iota}=\hat{\rho}\hat{\epsilon}\\
\Rightarrow\hat{\rho}=\lambda^{-1}(s+\lambda)\hat{\iota}\Rightarrow\rho=\lambda^{-1}\mathcal{L}^{-1}\left[s\hat{\iota}\right]+\iota\\
\Rightarrow\rho(x)=\frac{1}{\lambda}\left(\frac{d}{dx}\iota(x)+\iota(0)\right)+\iota(x),\label{eq:pprd_dead_time_match}\end{multline}
where $\hat{\,}$ decorates a function which was transformed by the
Laplace transform $\mathcal{L}$, and $s$ denotes the Laplace variable.
The renewal process defined by $\iota$ can be represented by a PPRD
if $\rho$ is a PDF. Let us call the hazard function of the renewal
process $h(x)$, and the survivor function $F(x)=\exp(-\int_{0}^{x}h(x^{\prime})dx^{\prime})$,
which obey $\iota(x)=h(x)F(x)$ \citep{Cox67}. Assume that $\iota(x)$
is differentiable. Since expression \eqref{eq:pprd_dead_time_match}
is always normalized, in order for it to define a suitable PDF we
only have to require $\rho(x)\geq0$ for all $x$, possibly in the
sense of distributions. This translates into\begin{equation}
\lambda^{-1}\left(h^{\prime}(x)-h^{2}(x)\right)+h(x)\geq0.\label{eq:pprd_hazard_condition}\end{equation}
In case $h(x)>0$, this can be written as \begin{equation}
h(x)-\frac{h^{\prime}(x)}{h(x)}\leq\lambda.\label{eq:pprd_hazard_condition_pos}\end{equation}
If, in addition, the hazard and its derivative are bounded in the
sense that $h(x)<\infty$ and $h^{\prime}(x)>-\infty$, there exists
a $\lambda>0$ such that \eqref{eq:pprd_hazard_condition_pos} is
fulfilled. These conditions are indeed met by a large class of renewal
processes.

For example, the gamma-process which has random inter-event intervals
with PDF $\iota(x)=\kappa_{r}(x)$ \eqref{eq:rho_gamma} with parameters
$r,\beta\in\mathbb{R}$, $r\geq1$, $\gammarate\geq0$ is by \eqref{eq:pprd_dead_time_match}
equivalent to the PPRD with  $\rho(x)=\kappa_{r-1}(x)$ and $\lambda=\beta$,
but other choices of $\lambda$ are also possible. This illustrates
the well known fact that the inter-event intervals of gamma processes
with integer shape parameter $n$ can be considered as the concatenation
of $n$ exponentially distributed intervals. In neuroscience the gamma-process
is frequently used to model stationary time series of action potential
emissions of nerve cells. To describe adaptation phenomena, a time-dependence
of the parameters of the hazard function was introduced in \citep{Muller07_2958}.
Identification of the gamma-process with a PPRD entails the alternative
to generalize the gamma process to time-dependent rates by varying
the input rate of the PPRD. Similarly, the log-normal process can
be represented as a PPRD. We define its inter-event interval as $x=\xi\Delta$,
where $\xi$ is a unit-less random number and $\Delta$ gives the
time-scale. Let $\xi$ be distributed according to the log-normal
PDF \begin{eqnarray*}
\eta(\xi) & = & \frac{1}{\sqrt{2\pi}\xi\sigma}\exp\left(-\frac{(\log\xi-\mu)^{2}}{2\sigma^{2}}\right)\end{eqnarray*}
for $\xi>0$, $\eta(0)=0$, with unit-less parameters $\mu,\,\sigma$.
Then $x$ is distributed according to $\iota(x)=\Delta^{-1}\eta(\nicefrac{x}{\Delta})$.
According to \eqref{eq:pprd_dead_time_match} the process can be represented
by any of the PPRDs with \begin{eqnarray*}
\rho(x) & = & \iota(x)\left(1-\frac{1}{x\lambda}\left(1+\frac{\log\frac{x}{\Delta}-\mu}{\sigma^{2}}\right)\right)\\
\lambda & \geq & \Delta^{-1}\sigma^{-2}\exp\left(-1-\mu+\sigma^{2}\right),\end{eqnarray*}
where the lower bound on $\lambda$ is due to the requirement $\rho(x)\geq0$.
For these and other renewal processes for which a PPRD representation
exists, non-equilibrium dynamics can be studied on the basis of \eqref{eq:DDE_pprd}.

\section{Discussion}

In this paper, we consider the effect of refractoriness on the output
of an encoding point process in case of arbitrary time-dependent input
signals. Such point processes, for example, are used to model the
generation of action potentials by nerve cells, the release and reuptake
of vesicles into the synaptic cleft, or the detection of particles
by technical devices. We describe ensembles of these stochastic processes
by the occupation numbers of two states: active and refractory. The
active components behave as inhomogeneous Poisson processes, but after
an event is produced the component is silent for the duration of the
dead-time, it is caught in a delay line. We derive a distributed delay
differential equation that describes the dynamics in the general case
of a randomly distributed dead-time.

Due to the simpler dynamics in case of a fixed dead-time, we first
elaborate properties of the PPD. For stationary input rate, we solve
the dynamics of the ensemble in a way that sheds light on the connection
between the fundamental solution of the DDE and the auto-correlation
function of the point process. This relation is employed to express
the time-dependent ensemble rate (output) for a step-change of the
hazard rate (input). The resulting output rate displays stochastic
transients and oscillations with a periodicity given by the dead-time.
Such transients might enable nerve cells to respond reliably to rapid
changes in the input currents \citep{Mainen95_1503,Berry98}. For
periodically modulated input rate, we demonstrate how the spectrum
of the steady-state periodic output rate results from the linear coupling
between harmonics. In the particular case of cosine-modulated input
signals only adjacent harmonics are coupled. This nearest-neighbor
interaction is rigorously solved using the theory of three-term-recurrence
relations and continued fractions \citep{Gautschi67}. 

Our analytic result explains frequency doubling, the emergence of
higher harmonics and the dependence of the time averaged population
activity on the modulation frequency. In particular, slow frequency
components of the input are attenuated and distorted in the population
rate, which is in contrast to the claim that the PPD transmits slow
frequency signals more reliably than the Poisson process \citep{Gerstner02}. 

In case of periodic input modulation, the output spectrum contains
all harmonics of the fundamental frequency of the input. This might
be related to a psychophysical phenomenon called {}``missing fundamental
illusion'' \citep{Schouten40,Licklider51}: Being presented an auditory
stimulus which consists of several harmonics of a fundamental frequency,
but in which the fundamental frequency itself is missing, subjects
nonetheless perceive the fundamental frequency as if it was contained
in the stimulus spectrum. By considering neurons in the auditory system
as PPDs whose hazard rate is modulated by the auditory stimulus, our
theory explains how the lowest harmonic is recovered in the population
activity of the neurons. Conversely, our results can be applied to
infer input rate profiles from the count rate of detectors with dead-time,
in particular in the case of periodic input, for which \eqref{eq:beta_k_system_per}
applies. 

For the more general case of a random, arbitrarily distributed dead-time,
we show how the DDE generalizes to a distributed DDE. By suitable
choice of the distribution of the dead-time, non-equilibrium dynamics
of a large class of renewal processes can be described. For integer
gamma-distributed dead-time we demonstrate how the distributed DDE
transforms into a coupled system of finitely many ordinary differential
equations, which could also be implemented as a multi-state Markov
system \citep{Toyoizumi09}. Regarding the output rate transient upon
a step change of the input and the transmission of periodic inputs,
we find that the qualitative behavior of the system is very similar
to the PPD. In conclusion, we present a canonical model for non-stationary
renewal processes, as well as the analytical methods to describe ensembles
thereof. 

\acknowledgments{We thank Eilif M\"uller for helpful discussions. Partially funded by BMBF Grant 01GQ0420 to BCCN Freiburg and EU Grant 15879 (FACETS).}

\bibliographystyle{../bib/apsrev4-1}
\bibliography{../bib/brain,../bib/computer,../bib/deger,../bib/math,../bib/BCCN}

\begin{thebibliography}{28}%
\makeatletter
\providecommand \@ifxundefined [1]{%
 \@ifx{#1\undefined}
}%
\providecommand \@ifnum [1]{%
 \ifnum #1\expandafter \@firstoftwo
 \else \expandafter \@secondoftwo
 \fi
}%
\providecommand \@ifx [1]{%
 \ifx #1\expandafter \@firstoftwo
 \else \expandafter \@secondoftwo
 \fi
}%
\providecommand \natexlab [1]{#1}%
\providecommand \enquote  [1]{``#1''}%
\providecommand \bibnamefont  [1]{#1}%
\providecommand \bibfnamefont [1]{#1}%
\providecommand \citenamefont [1]{#1}%
\providecommand \href@noop [0]{\@secondoftwo}%
\providecommand \href [0]{\begingroup \@sanitize@url \@href}%
\providecommand \@href[1]{\@@startlink{#1}\@@href}%
\providecommand \@@href[1]{\endgroup#1\@@endlink}%
\providecommand \@sanitize@url [0]{\catcode `\\12\catcode `\$12\catcode
  `\&12\catcode `\#12\catcode `\^12\catcode `\_12\catcode `\%12\relax}%
\providecommand \@@startlink[1]{}%
\providecommand \@@endlink[0]{}%
\providecommand \url  [0]{\begingroup\@sanitize@url \@url }%
\providecommand \@url [1]{\endgroup\@href {#1}{\urlprefix }}%
\providecommand \urlprefix  [0]{URL }%
\providecommand \Eprint [0]{\href }%
\@ifxundefined \urlstyle {%
  \providecommand \doi  [0]{\begingroup \@sanitize@url \@doi}%
  \providecommand \@doi [1]{\endgroup \@@startlink {\doibase
  #1}doi:\discretionary {}{}{}#1\@@endlink }%
}{%
  \providecommand \doi  [0]{doi:\discretionary{}{}{}\begingroup
  \urlstyle{rm}\Url }%
}%
\providecommand \doibase [0]{http://dx.doi.org/}%
\providecommand \Doi [0]{\begingroup \@sanitize@url \@Doi }%
\providecommand \@Doi  [1]{\endgroup\@@startlink{\doibase#1}\@@Doi}%
\providecommand \@@Doi [1]{#1\@@endlink}%
\providecommand \selectlanguage [0]{\@gobble}%
\providecommand \bibinfo  [0]{\@secondoftwo}%
\providecommand \bibfield  [0]{\@secondoftwo}%
\providecommand \translation [1]{[#1]}%
\providecommand \BibitemOpen [0]{}%
\providecommand \bibitemStop [0]{}%
\providecommand \bibitemNoStop [0]{.\EOS\space}%
\providecommand \EOS [0]{\spacefactor3000\relax}%
\providecommand \BibitemShut  [1]{\csname bibitem#1\endcsname}%
\bibitem [{\citenamefont {Cox}(1967)}]{Cox67}%
  \BibitemOpen
  \bibfield  {author} {\bibinfo {author} {\bibfnamefont {D.}~\bibnamefont
  {Cox}},\ }\href@noop {} {\emph {\bibinfo {title} {Renewal Theory}}}\
  (\bibinfo  {publisher} {London: Chapman and Hall},\ \bibinfo {address}
  {Science Paperbacks},\ \bibinfo {year} {1967})\BibitemShut {NoStop}%
\bibitem [{\citenamefont {Lowen}\ and\ \citenamefont
  {Teich}(2005)}]{LowenTeich05}%
  \BibitemOpen
  \bibfield  {author} {\bibinfo {author} {\bibfnamefont {S.~B.}\ \bibnamefont
  {Lowen}}\ and\ \bibinfo {author} {\bibfnamefont {M.~C.}\ \bibnamefont
  {Teich}},\ }\href@noop {} {\emph {\bibinfo {title} {Fractal-Based Point
  Processes}}}\ (\bibinfo  {publisher} {John Wiley \& Sons},\ \bibinfo {year}
  {2005})\BibitemShut {NoStop}%
\bibitem [{\citenamefont {Grupen}\ and\ \citenamefont
  {Shwartz}(2008)}]{Grupen08}%
  \BibitemOpen
  \bibfield  {author} {\bibinfo {author} {\bibfnamefont {C.}~\bibnamefont
  {Grupen}}\ and\ \bibinfo {author} {\bibfnamefont {B.}~\bibnamefont
  {Shwartz}},\ }\href@noop {} {\emph {\bibinfo {title} {Particle detectors}}}\
  (\bibinfo  {publisher} {Cambridge University Press},\ \bibinfo {year}
  {2008})\ ISBN \bibinfo {isbn} {978-0-521-84006-4}\BibitemShut {NoStop}%
\bibitem [{\citenamefont {M\"uller}(1981)}]{Muller81}%
  \BibitemOpen
  \bibinfo {editor} {\bibfnamefont {J.~W.}\ \bibnamefont {M\"uller}},\ ed.,\
  \href@noop {} {\emph {\bibinfo {title} {Bibliography on dead time effects.
  Report BIPM-81/11}}}\ (\bibinfo  {publisher} {Bureau International des Poids
  et Measures, S\`evres, France},\ \bibinfo {year} {1981})\BibitemShut
  {NoStop}%
\bibitem [{\citenamefont {Loebel}\ \emph {et~al.}(2009)\citenamefont {Loebel},
  \citenamefont {Silberberg}, \citenamefont {Helbig}, \citenamefont {Markram},
  \citenamefont {Tsodyks},\ and\ \citenamefont {Richardson}}]{Loebel09}%
  \BibitemOpen
  \bibfield  {author} {\bibinfo {author} {\bibfnamefont {A.}~\bibnamefont
  {Loebel}}, \bibinfo {author} {\bibfnamefont {G.}~\bibnamefont {Silberberg}},
  \bibinfo {author} {\bibfnamefont {D.}~\bibnamefont {Helbig}}, \bibinfo
  {author} {\bibfnamefont {H.}~\bibnamefont {Markram}}, \bibinfo {author}
  {\bibfnamefont {M.}~\bibnamefont {Tsodyks}}, \ and\ \bibinfo {author}
  {\bibfnamefont {M.~J.~E.}\ \bibnamefont {Richardson}},\ }\href@noop {}
  {\bibfield  {journal} {\bibinfo  {journal} {Front. Comput. Neurosci.},\
  }\textbf {\bibinfo {volume} {3}} (\bibinfo {year} {2009})}\BibitemShut
  {NoStop}%
\bibitem [{\citenamefont {Gerstner}\ and\ \citenamefont
  {Kistler}(2002)}]{Gerstner02}%
  \BibitemOpen
  \bibfield  {author} {\bibinfo {author} {\bibfnamefont {W.}~\bibnamefont
  {Gerstner}}\ and\ \bibinfo {author} {\bibfnamefont {W.}~\bibnamefont
  {Kistler}},\ }\href@noop {} {\emph {\bibinfo {title} {Spiking Neuron Models:
  Single Neurons, Populations, Plasticity}}}\ (\bibinfo  {publisher} {Cambridge
  University Press},\ \bibinfo {year} {2002})\ ISBN \bibinfo {isbn}
  {(paperback) 0521890799}\BibitemShut {NoStop}%
\bibitem [{\citenamefont {Kuhn}\ \emph {et~al.}(2004)\citenamefont {Kuhn},
  \citenamefont {Aertsen},\ and\ \citenamefont {Rotter}}]{Kuhn04}%
  \BibitemOpen
  \bibfield  {author} {\bibinfo {author} {\bibfnamefont {A.}~\bibnamefont
  {Kuhn}}, \bibinfo {author} {\bibfnamefont {A.}~\bibnamefont {Aertsen}}, \
  and\ \bibinfo {author} {\bibfnamefont {S.}~\bibnamefont {Rotter}},\
  }\href@noop {} {\bibfield  {journal} {\bibinfo  {journal} {J. Neurosci.},\
  }\textbf {\bibinfo {volume} {24}},\ \bibinfo {pages} {2345} (\bibinfo {year}
  {2004})}\BibitemShut {NoStop}%
\bibitem [{\citenamefont {Brenguier}\ and\ \citenamefont
  {Amodei}(1989)}]{Brenguier89}%
  \BibitemOpen
  \bibfield  {author} {\bibinfo {author} {\bibfnamefont {J.}~\bibnamefont
  {Brenguier}}\ and\ \bibinfo {author} {\bibfnamefont {L.}~\bibnamefont
  {Amodei}},\ }\Doi {10.1175/1520-0426(1989)006<0575:CADTCF>2.0.CO;2}
  {\bibfield  {journal} {\bibinfo  {journal} {J. Atmos. Ocean. Tech.},\
  }\textbf {\bibinfo {volume} {6}},\ \bibinfo {pages} {575} (\bibinfo {year}
  {1989})}\BibitemShut {NoStop}%
\bibitem [{\citenamefont {Yu}\ and\ \citenamefont {Fessler}(2000)}]{Yu00}%
  \BibitemOpen
  \bibfield  {author} {\bibinfo {author} {\bibfnamefont {D.~F.}\ \bibnamefont
  {Yu}}\ and\ \bibinfo {author} {\bibfnamefont {J.~A.}\ \bibnamefont
  {Fessler}},\ }\Doi {10.1016/S0168-9002(02)00460-6} {\bibfield  {journal}
  {\bibinfo  {journal} {Phys. Med. Biol.},\ }\textbf {\bibinfo {volume} {45}},\
  \bibinfo {pages} {2043} (\bibinfo {year} {2000})}\BibitemShut {NoStop}%
\bibitem [{\citenamefont {Teich}\ and\ \citenamefont {McGill}(1976)}]{Teich76}%
  \BibitemOpen
  \bibfield  {author} {\bibinfo {author} {\bibfnamefont {M.~C.}\ \bibnamefont
  {Teich}}\ and\ \bibinfo {author} {\bibfnamefont {W.~J.}\ \bibnamefont
  {McGill}},\ }\href@noop {} {\bibfield  {journal} {\bibinfo  {journal} {Phys.
  Rev. Lett.},\ }\textbf {\bibinfo {volume} {36}},\ \bibinfo {pages} {754}
  (\bibinfo {year} {1976})}\BibitemShut {NoStop}%
\bibitem [{\citenamefont {Picinbono}(2009)}]{Picinbono09}%
  \BibitemOpen
  \bibfield  {author} {\bibinfo {author} {\bibfnamefont {B.}~\bibnamefont
  {Picinbono}},\ }\href@noop {} {\bibfield  {journal} {\bibinfo  {journal}
  {Commun. Stat. Simulat.},\ }\textbf {\bibinfo {volume} {38}},\ \bibinfo
  {pages} {2198} (\bibinfo {year} {2009})}\BibitemShut {NoStop}%
\bibitem [{\citenamefont {Vannucci}\ and\ \citenamefont
  {Teich}(1978)}]{Vanucci78}%
  \BibitemOpen
  \bibfield  {author} {\bibinfo {author} {\bibfnamefont {G.}~\bibnamefont
  {Vannucci}}\ and\ \bibinfo {author} {\bibfnamefont {M.~C.}\ \bibnamefont
  {Teich}},\ }\href@noop {} {\bibfield  {journal} {\bibinfo  {journal} {Opt.
  Commun.},\ }\textbf {\bibinfo {volume} {25}},\ \bibinfo {pages} {267}
  (\bibinfo {year} {1978})}\BibitemShut {NoStop}%
\bibitem [{\citenamefont {Larsen}\ and\ \citenamefont
  {Kostinski}(2009)}]{Larsen09}%
  \BibitemOpen
  \bibfield  {author} {\bibinfo {author} {\bibfnamefont {M.~L.}\ \bibnamefont
  {Larsen}}\ and\ \bibinfo {author} {\bibfnamefont {A.~B.}\ \bibnamefont
  {Kostinski}},\ }\Doi {10.1088/0957-0233/20/9/095101} {\bibfield  {journal}
  {\bibinfo  {journal} {Meas. Sci. Technol.},\ }\textbf {\bibinfo {volume}
  {20}} (\bibinfo {year} {2009})},\ \doi
  {10.1088/0957-0233/20/9/095101}\BibitemShut {NoStop}%
\bibitem [{\citenamefont {Turcott}\ \emph {et~al.}(1994)\citenamefont
  {Turcott}, \citenamefont {Lowen}, \citenamefont {Li}, \citenamefont
  {Johnson}, \citenamefont {Tsuchitani},\ and\ \citenamefont
  {Teich}}]{Turcott94}%
  \BibitemOpen
  \bibfield  {author} {\bibinfo {author} {\bibfnamefont {R.~G.}\ \bibnamefont
  {Turcott}}, \bibinfo {author} {\bibfnamefont {S.~B.}\ \bibnamefont {Lowen}},
  \bibinfo {author} {\bibfnamefont {E.}~\bibnamefont {Li}}, \bibinfo {author}
  {\bibfnamefont {D.~H.}\ \bibnamefont {Johnson}}, \bibinfo {author}
  {\bibfnamefont {C.}~\bibnamefont {Tsuchitani}}, \ and\ \bibinfo {author}
  {\bibfnamefont {M.}~\bibnamefont {Teich}},\ }\href@noop {} {\bibfield
  {journal} {\bibinfo  {journal} {Biol. Cybern.},\ }\textbf {\bibinfo {volume}
  {70}},\ \bibinfo {pages} {209} (\bibinfo {year} {1994})}\BibitemShut
  {NoStop}%
\bibitem [{\citenamefont {Berry}\ and\ \citenamefont
  {Meister}(1998)}]{Berry98}%
  \BibitemOpen
  \bibfield  {author} {\bibinfo {author} {\bibfnamefont {M.~J.}\ \bibnamefont
  {Berry}}\ and\ \bibinfo {author} {\bibfnamefont {M.}~\bibnamefont
  {Meister}},\ }\href@noop {} {\bibfield  {journal} {\bibinfo  {journal} {J.
  Neurosci.},\ }\textbf {\bibinfo {volume} {18}},\ \bibinfo {pages} {2200}
  (\bibinfo {year} {1998})}\BibitemShut {NoStop}%
\bibitem [{\citenamefont {Johnson}\ and\ \citenamefont
  {Swami}(1983)}]{Johnson83}%
  \BibitemOpen
  \bibfield  {author} {\bibinfo {author} {\bibfnamefont {D.~H.}\ \bibnamefont
  {Johnson}}\ and\ \bibinfo {author} {\bibfnamefont {A.}~\bibnamefont
  {Swami}},\ }\href@noop {} {\bibfield  {journal} {\bibinfo  {journal} {J.
  Acoust. Soc. America},\ }\textbf {\bibinfo {volume} {74}},\ \bibinfo {pages}
  {493} (\bibinfo {year} {1983})}\BibitemShut {NoStop}%
\bibitem [{\citenamefont {Wilson}\ and\ \citenamefont
  {Cowan}(1972)}]{Wilson72_1}%
  \BibitemOpen
  \bibfield  {author} {\bibinfo {author} {\bibfnamefont {H.~R.}\ \bibnamefont
  {Wilson}}\ and\ \bibinfo {author} {\bibfnamefont {J.~D.}\ \bibnamefont
  {Cowan}},\ }\href@noop {} {\bibfield  {journal} {\bibinfo  {journal}
  {Biophys. J.},\ }\textbf {\bibinfo {volume} {12}},\ \bibinfo {pages} {1}
  (\bibinfo {year} {1972})}\BibitemShut {NoStop}%
\bibitem [{\citenamefont {Brown}\ \emph {et~al.}(2001)\citenamefont {Brown},
  \citenamefont {Barbieri}, \citenamefont {Ventura}, \citenamefont {Kaas},\
  and\ \citenamefont {Frank}}]{Brown01}%
  \BibitemOpen
  \bibfield  {author} {\bibinfo {author} {\bibfnamefont {E.~N.}\ \bibnamefont
  {Brown}}, \bibinfo {author} {\bibfnamefont {R.}~\bibnamefont {Barbieri}},
  \bibinfo {author} {\bibfnamefont {V.}~\bibnamefont {Ventura}}, \bibinfo
  {author} {\bibfnamefont {R.~E.}\ \bibnamefont {Kaas}}, \ and\ \bibinfo
  {author} {\bibfnamefont {L.~M.}\ \bibnamefont {Frank}},\ }\href@noop {}
  {\bibfield  {journal} {\bibinfo  {journal} {Neural Comput.},\ }\textbf
  {\bibinfo {volume} {14}},\ \bibinfo {pages} {325} (\bibinfo {year}
  {2001})}\BibitemShut {NoStop}%
\bibitem [{\citenamefont {Reich}\ \emph {et~al.}(1998)\citenamefont {Reich},
  \citenamefont {Victor},\ and\ \citenamefont {Knight}}]{Reich98}%
  \BibitemOpen
  \bibfield  {author} {\bibinfo {author} {\bibfnamefont {D.~S.}\ \bibnamefont
  {Reich}}, \bibinfo {author} {\bibfnamefont {J.~D.}\ \bibnamefont {Victor}}, \
  and\ \bibinfo {author} {\bibfnamefont {B.~W.}\ \bibnamefont {Knight}},\
  }\href@noop {} {\bibfield  {journal} {\bibinfo  {journal} {J. Neurosci.},\
  }\textbf {\bibinfo {volume} {18}},\ \bibinfo {pages} {10090} (\bibinfo {year}
  {1998})}\BibitemShut {NoStop}%
\bibitem [{\citenamefont {Nawrot}\ \emph {et~al.}(2008)\citenamefont {Nawrot},
  \citenamefont {Boucsein}, \citenamefont {Rodriguez~Molina}, \citenamefont
  {Riehle}, \citenamefont {Aertsen},\ and\ \citenamefont
  {Rotter}}]{Nawrot08_374}%
  \BibitemOpen
  \bibfield  {author} {\bibinfo {author} {\bibfnamefont {M.~P.}\ \bibnamefont
  {Nawrot}}, \bibinfo {author} {\bibfnamefont {C.}~\bibnamefont {Boucsein}},
  \bibinfo {author} {\bibfnamefont {V.}~\bibnamefont {Rodriguez~Molina}},
  \bibinfo {author} {\bibfnamefont {A.}~\bibnamefont {Riehle}}, \bibinfo
  {author} {\bibfnamefont {A.}~\bibnamefont {Aertsen}}, \ and\ \bibinfo
  {author} {\bibfnamefont {S.}~\bibnamefont {Rotter}},\ }\href@noop {}
  {\bibfield  {journal} {\bibinfo  {journal} {J. Neurosci. Methods},\ }\textbf
  {\bibinfo {volume} {169}},\ \bibinfo {pages} {374} (\bibinfo {year}
  {2008})}\BibitemShut {NoStop}%
\bibitem [{\citenamefont {Muller}\ \emph {et~al.}(2007)\citenamefont {Muller},
  \citenamefont {Buesing}, \citenamefont {Schemmel},\ and\ \citenamefont
  {Meier}}]{Muller07_2958}%
  \BibitemOpen
  \bibfield  {author} {\bibinfo {author} {\bibfnamefont {E.}~\bibnamefont
  {Muller}}, \bibinfo {author} {\bibfnamefont {L.}~\bibnamefont {Buesing}},
  \bibinfo {author} {\bibfnamefont {J.}~\bibnamefont {Schemmel}}, \ and\
  \bibinfo {author} {\bibfnamefont {K.}~\bibnamefont {Meier}},\ }\href@noop {}
  {\bibfield  {journal} {\bibinfo  {journal} {Neural Comput.},\ }\textbf
  {\bibinfo {volume} {19}},\ \bibinfo {pages} {2958} (\bibinfo {year}
  {2007})}\BibitemShut {NoStop}%
\bibitem [{\citenamefont {Toyoizumi}\ \emph {et~al.}(2009)\citenamefont
  {Toyoizumi}, \citenamefont {Rad},\ and\ \citenamefont
  {Paninski}}]{Toyoizumi09}%
  \BibitemOpen
  \bibfield  {author} {\bibinfo {author} {\bibfnamefont {T.}~\bibnamefont
  {Toyoizumi}}, \bibinfo {author} {\bibfnamefont {K.~R.}\ \bibnamefont {Rad}},
  \ and\ \bibinfo {author} {\bibfnamefont {L.}~\bibnamefont {Paninski}},\
  }\href@noop {} {\bibfield  {journal} {\bibinfo  {journal} {Neural Comput.},\
  }\textbf {\bibinfo {volume} {21}},\ \bibinfo {pages} {1203} (\bibinfo {year}
  {2009})}\BibitemShut {NoStop}%
\bibitem [{\citenamefont {Bellman}\ and\ \citenamefont
  {Cooke}(1963)}]{Bellman63}%
  \BibitemOpen
  \bibfield  {author} {\bibinfo {author} {\bibfnamefont {R.~E.}\ \bibnamefont
  {Bellman}}\ and\ \bibinfo {author} {\bibfnamefont {K.~L.}\ \bibnamefont
  {Cooke}},\ }\href@noop {} {\emph {\bibinfo {title} {Differential-Difference
  Equations}}}\ (\bibinfo  {publisher} {RAND Corporation},\ \bibinfo {year}
  {1963})\BibitemShut {NoStop}%
\bibitem [{\citenamefont {Gautschi}(1967)}]{Gautschi67}%
  \BibitemOpen
  \bibfield  {author} {\bibinfo {author} {\bibfnamefont {W.}~\bibnamefont
  {Gautschi}},\ }\href@noop {} {\bibfield  {journal} {\bibinfo  {journal} {SIAM
  Rev.},\ }\textbf {\bibinfo {volume} {9}},\ \bibinfo {pages} {24} (\bibinfo
  {year} {1967})}\BibitemShut {NoStop}%
\bibitem [{\citenamefont {Ripley}(1987)}]{Ripley87}%
  \BibitemOpen
  \bibfield  {author} {\bibinfo {author} {\bibfnamefont {B.~D.}\ \bibnamefont
  {Ripley}},\ }\href@noop {} {\emph {\bibinfo {title} {Stochastic
  Simulation}}},\ Wiley Series in Probability and Mathematical Statistics\
  (\bibinfo  {publisher} {John Wiley \& Sons},\ \bibinfo {address} {New York},\
  \bibinfo {year} {1987})\BibitemShut {NoStop}%
\bibitem [{\citenamefont {Mainen}\ and\ \citenamefont
  {Sejnowski}(1995)}]{Mainen95_1503}%
  \BibitemOpen
  \bibfield  {author} {\bibinfo {author} {\bibfnamefont {Z.~F.}\ \bibnamefont
  {Mainen}}\ and\ \bibinfo {author} {\bibfnamefont {T.~J.}\ \bibnamefont
  {Sejnowski}},\ }\href@noop {} {\bibfield  {journal} {\bibinfo  {journal}
  {Science},\ }\textbf {\bibinfo {volume} {268}},\ \bibinfo {pages} {1503}
  (\bibinfo {year} {1995})}\BibitemShut {NoStop}%
\bibitem [{\citenamefont {Schouten}(1940)}]{Schouten40}%
  \BibitemOpen
  \bibfield  {author} {\bibinfo {author} {\bibfnamefont {J.}~\bibnamefont
  {Schouten}},\ }\href@noop {} {\bibfield  {journal} {\bibinfo  {journal}
  {Philips tech. Rev. 5}} (\bibinfo {year} {1940})}\BibitemShut {NoStop}%
\bibitem [{\citenamefont {Licklider}(1951)}]{Licklider51}%
  \BibitemOpen
  \bibfield  {author} {\bibinfo {author} {\bibfnamefont {J.~C.~R.}\
  \bibnamefont {Licklider}},\ }\href@noop {} {\bibfield  {journal} {\bibinfo
  {journal} {Cell. Mol. Life Sci.},\ }\textbf {\bibinfo {volume} {7}},\
  \bibinfo {pages} {128} (\bibinfo {year} {1951})}\BibitemShut {NoStop}%
\end{thebibliography}%

\end{document}